\documentclass[12pt]{article}
\usepackage{amsmath}
\usepackage{amssymb}
\usepackage{amsthm}
\usepackage{amsfonts}
\usepackage{graphicx}

\usepackage{bbm,epsfig,graphics,epic,color,rotating,color}
\numberwithin{equation}{section}

\newcommand{\R}{{\mathbb R}}
\newcommand{\Z}{{\mathbb Z}}

\newcommand{\un}[1]{\underline{#1}}

\newcommand{\ed}{\mathrm{d}}

\newcommand{\reff}[1]{(\ref{#1})}

\newtheorem{theorem}{Theorem}[section]
\newtheorem{lemma}[theorem]{Lemma}

\newtheorem{remark}[theorem]{Remark}

\theoremstyle{definition}

\title{
Large deviations for excursions of non-homogeneous Markov processes}

\author{ A. Mogulskii$^{1}$ \and E. Pechersky~$^{2}$\and A. Yambartsev~$^{3}$
}
\date{}

\begin{document}

\maketitle {\footnotesize
\noindent $^{1}$ Sobolev Institute of Mathematics of  Siberian Branch of Russian Academy of Sciences,\\
4, ac. Koptyug str., Novosibirsk, Russia,\\
E-mail: mogul@math.nsc.ru

\noindent $^2$ Dobrushin laboratory of Institute for Information
Transmission Problems of Russian Academy of Sciences,\\
19, Bolshoj Karetny, Moscow, Russia.\\
E-mail: pech@iitp.ru

\noindent $^3$ Department of Statistics, Institute of Mathematics
and Statistics, University of S\~ao Paulo, Rua do Mat\~ao 1010,
CEP 05508--090, S\~ao Paulo SP, Brazil.\\
E-mail: yambar@ime.usp.br }

\begin{abstract}
In this paper, the large deviations on trajectory level for ergodic Markov processes are studied. These processes take values in the non-negative quadrant of the two dimension lattice and are concentrated on step-wise functions. The rates of jumps towards the axes (jump down) depend on the position of the process  -- the higher the position, the greater the rate. The rates of jumps going in the same direction as the axes (jump up) are constants. Therefore the processes are ergodic. The large deviations are studied under equal scalings of both space and time. The scaled versions of the processes converge to 0. The main result is that the probabilities of long excursions out of 0 tend to 0 exponentially fast with an exponent proportional to the square of the scaling parameter. A proportionality coefficient is an integral of a linear combination of path components. A rate function of the large deviation principle is calculated for continuous functions only.
\end{abstract}

\section{Introduction} There are different settings in the large deviation theory studying probabilities of rare events (see, for example, the books \cite{DZ,DS,H, OV,Ph,V}).

This paper is devoted to investigations of the rare event probabilities for a class of ergodic Markov processes. The goal is to find an asymptotic behavior of logarithm of probabilities for a long excursion of the process far from equilibrium states. 
We apply the large deviation setting using equal contractions in time and in space. The path level of a large deviation principle is obtained.

A basic random object is a continuous time Markov ergodic process $\xi$ with state space
$\Z^2_+:=\{(z_1,z_2):~z_1\ge 0,~z_2\ge
0\}$. Paths of
$\xi$  are  piece-wise
constant functions.
The jumps belong to the following set
$$
\mathcal{Y}=\{(1,0),(0,1),(-1,0),(0,-1),(-1,-1)\}.
$$
The probabilities of the jumps are such that they do not take the process outside of $\Z^2_+$. The intensities of
the jumps depend on the value of $\xi$ at the moment before the jump.
If at a   moment $t$  the process  value is equal to
$\xi(t)=(z_1,z_2)$, then any enlargement (jump up) of at least one of the components of $(z_1,z_2)$ happens
with a constant intensity. However, any reduction (jump down) of at least one of the components of $(z_1,z_2)$
happens with an intensity proportional to this  co-ordinate. This
property implies the ergodicity of the process $\xi$.

We consider the large deviations for the sequence
$\xi_T(t)=\left(\frac{\xi(tT)}{T}\right)_{T>0}$ of the processes on
$t\in[0,1]$ with $\xi(0)=(0,0)$. The large deviation principle for $\xi$ is being established on a  set of c\`{a}dl\`{a}g functions $X$ with a finite number of their jumps, which includes all typical paths of $\xi$. A rate function is finite for a set $F$ of continuous functions on $[0,1]$ such that any $\un{f}\in F$ has positive co-ordinates except of their values at $t=0$, where $\un{f}(0)=(0,0)$. When the processes $\xi_T$ are localized in a small neighborhood of some function $\un f\in F$, we say that the process $\xi$ has a long excursion far from the equilibrium. We find that the rate function of $\un f=(f_1,f_2)$ is the following integral form
\begin{equation}\label{rateI}
I(\un{f})=\int_0^1 \bigl( c_1f_1(t)+c_2f_2(t)+c_3\min\{f_1(t),f_2(t)\} \bigr) dt,
\end{equation}
where constants $c_1,c_2$ and $c_3$ are parameters defining the process $\xi$ (see exact definitions in section \ref{22}). A local principle of the large deviations proved in this paper implies that the probability of a long excursion in a small neighborhood $U(\un f)$ of a function $\un f\in F$ has an order
\begin{equation*}
    e^{-T^2I(\un{f})}.
\end{equation*}
We use in this paper the uniform topology in $F$.

Notice that  derivatives of $\un f$ are not included in the expression for $I(\un f)$ \reff{rateI}.
Such form of the rate function seems paradoxical. Indeed, let a continuous function $g_1:\:[0,1]\to\R_+$ have a form of a high narrow peak such that $\int_0^1g_1(t)\ed t=\varepsilon_0$ is small, and let $\un{g}=(g_1,0)$. The difference of the rate functions $I(\un{f}+\un{g})$ and $I(\un{f})$, for $\un f \in F$, is small and equal to $c_1\varepsilon_0$, but $\sup_{t} \{g_1(t)-f_1(t)\}$ can be very large.
An explication of this paradox is that the probability that the process $\xi_T$ belongs to a ``neighborhood" of $\un{g}$ has an order
\begin{equation}\label{stac}
e^{-T\ln(T)C},
\end{equation}
where $C$ is a constant which dependents on $\un{g}$.
The asymptotic \reff{stac} is not proved in this paper. The word ``neighborhood" is under quotes because \reff{stac} has to be proved in  different settings (it will be done in another paper).
It shows that the probability of $\xi_T$ being out of zero for a long time is much less than the probability of a high ejection during a short period. 

This study was inspired by the work \cite{BSV}, where ergodic
properties of more complicated processes were studied. The goal of
the authors of \cite{BSV} was to describe market dynamics.
Our goal is focused on some peculiarities of the large deviations
for similar models and our  version of the model  is hardly  proper for  market investigations.

\section{Results.}

\subsection{Notations.}
Let $\xi(t)=(\xi_1(t),~\xi_2(t)), \ t\in [0,\infty)$ be a Markov process
with state space ${\mathbb Z}^2_+:=\{(z_1,z_2):~z_1\ge 0,~z_2\ge
0\}$. The evolution of the process can be described in the following
way. Let a state of the process at a moment $t\ge 0$ be
$\xi(t)=\un{z}=(z_1,z_2)\in\Z^2_+$. The state is not changed during a
time $\tau_{\un{z}}$, where   $\tau_{\un{z}}$ is a random variable
distributed exponentially with a parameter $h(\un{z})$. At the moment $t+\tau_{\un{z}}$
the value of the process becomes equal to
$\un{z}+\un{y}$, where $\un{y}$ belongs to
\begin{equation}\label{jump}
\mathcal{Y}= \{(1,0),(0,1),(-1,0),(0,-1),(-1,-1)\}.
\end{equation}
The intensities of the jumps is a sum
\begin{equation}\label{jum}
 h(\un{z}):= \lambda_{\un{z}}(1,0)+\lambda_{\un{z}}(0,1)+\lambda_{\un{z}}(-1,0)+\lambda_{\un{z}}(0,-1)+
  \lambda_{\un{z}}(-1,-1),
\end{equation}
where
\begin{eqnarray}
&&\lambda_{\un{z}}(1,0):=\lambda(1,0),~~\lambda_{\un{z}}(0,1):=\lambda(0,1),~~\nonumber\\
&&\lambda_{\un{z}}(-1,0):=z_1\lambda(-1,0),~~\lambda_{\un{z}}(0,-1):=z_2\lambda(0,-1),~~\label{int1}\\
&&\lambda_{\un{z}}(-1,-1):=\min\{z_1,z_2\}\lambda(-1,-1),\nonumber
\end{eqnarray}
and the constants $\lambda(\un y)$ at $\un y\in {\mathcal Y}$ are positive. The probability of the jump $\un{y}$ is
\begin{equation}\label{prj}
  p_{\un{z}}(\un{y}):=\frac{\lambda_{\un{z}}(\un{y})}{h(\un{z})},~~~\un{y}=(y_1,y_2)\in \mathcal{Y}.
\end{equation}

\subsection{The local large deviation  principle.}\label{22}
In this section we study the local deviation principle for the measures $(P_T)$ which are the
distributions of the processes $(\xi_T(t)=\frac{1}{T}\xi(tT)),\ t\in [0, 1]$. The support of the processes $\xi_T$ is a subset of the set $X$ of non-negative  c\`{a}dl\`{a}g functions
\begin{equation*}
\un x:\:[0,1]\to \R^2_+ = \{(y_1,y_2)\in \mathbb R^2: \ y_1 \ge 0, y_2 \ge 0\},
\end{equation*}
which are right-continuous and have left limits everywhere, having finite number of jumps on $[0,1]$ and such that $\un x (0) = (0,0)$ (definition of the c\`{a}dl\`{a}g functions see, for example, in \cite{B}).
We introduce an \textit{uniform} topology on $X$, which, in this case, is determined by the distance $d(\un x_1,\un x_2)$ between two functions $\un x_1,\un x_2\in X$ as follows
\begin{equation}\label{metr}
d(\un x_1,\un x_2)=\sup_{t\in[0,1]}\|\un x_1(t)-\un x_2(t)\|,
\end{equation}
where $\|\cdot\|$ means the usual Euclidean norm in $\R^2$.

There is a weak convergence $P_T\Rightarrow\delta_{\un x_0}$,
where $\un x_0(t)\equiv 0,~t\in[0,1]$.
Studying the long excursion far from $\un x_0$ we consider the set $F\subset X$ of continuous functions $\un{f}(t) = (f_1(t),f_2(t))$ satisfying the following properties:
\begin{description}
\item [$F_1$] $\un{f}(0)=(0,0)$,
\item [$F_2$] $f_1(t)>0$ and $f_2(t)>0$ for any $t > 0$.
\end{description}
We have found the rate function for this class $F\subset X$ of continuous functions satisfying the conditions $F_1$ and $F_2$.

For brevity we shall use  the notations $c_0=\lambda(1,0)+\lambda(0,1),\ c_1=\lambda(-1,0), \ c_2=\lambda(0,-1), \ c_3=\lambda(-1,-1)$.
Thus we rewrite \reff{jum} as (see also \reff{int1})
\begin{equation}\label{jum1}
h(\un{z}) \equiv h{(z_1,z_2)} = c_0+c_1z_1+c_2z_2+c_3\min\{z_1,z_2\}.
\end{equation}
On the set $X$ we define the following functional $I:\ X\to \mathbb R \cup \{\infty\}$
\begin{equation}\label{rate}
I(\un{x}):= \begin{cases} \int_0^1 \bigl( c_1x_1(t)+c_2x_2(t)+c_3\min\{x_1(t),x_2(t)\} \bigr) dt, & \mbox{ if } \un x \in F, \\
\infty, & \mbox{ if } \un x \notin F. \end{cases} 
\end{equation}
$I(\un{x})$ is finite for all bounded continuous functions $\un{x} \in F$. In the next theorem we prove the local large deviation principle with rate function $I(\un x)$.
\begin{theorem}\label{thlldp} For any $\un{f}\in F$
\begin{equation}\label{1}
 \lim_{\varepsilon\to 0} \lim_{T\to \infty}\frac{1}{T^2}\ln
  {\bf P} \bigl( \xi_T\in U_{\varepsilon}(\un{f}) \bigr) =-I(\un{f}),
\end{equation}
where (see \reff{metr})
\begin{equation}\label{ne}
U_\varepsilon(\un f)=\{\un g \in X:\: d(\un{f}, \un{g})<\varepsilon\}.
\end{equation}
\end{theorem}
\proof  \underline{Upper bound}. We have to show that
\begin{equation}\label{2}
L_+:=\limsup_{\varepsilon\to 0} \limsup_{T\to \infty}\frac{1}{T^2}\ln
  {\bf P} \bigl( \xi_T\in U_{\varepsilon}(\un{f}) \bigr) \le -I(\un{f}).
\end{equation}
In order to show it, consider a Markov process $\zeta(t)=(\zeta_1(t),\zeta_2(t)),~t\in [0,T]$, with 
state space ${\mathbb Z}^2$ and its intensity of jumps equal to 1. The
process $\zeta(t)$ is homogenous in time. At a jump moment the process $\zeta$ changes its value from $\un{z} \in \Z^2$ to $\un{z}+\un{y}$
with uniform probabilities $1/5$ for $\un{y}\in \mathcal Y$. 
It means that the process $\zeta$ is homogeneous in space, as well. The process $\zeta$ may be out of $\Z^2_ +$, moreover the process leaves $\Z_+^2$ with probability 1. 

Let $X_T$ be the set of all trajectories of the process $\xi$ on the time interval $[0,T]$.
The distribution of the
process $\xi$ is absolutely continuous with respect to $\zeta$ with  density
\begin{eqnarray}\nonumber
 {\cal P} \left( \un{u}(\cdot) \right) &=& 5^{N_T(\un{u})}
 \prod_{i=0}^{N_T(\un{u})-1} h(\un{u}(t_i))e^{-(h(\un{u}(t_i))-1)\tau_{i+1}} p_{ \un{u}(t_{i}) } (\un{u}(t_{i+1}) - \un{u}(t_i))\times \nonumber\\
&&h(\un u(t_{N_T(\un u)}))e^{-(h(\un u(t_{N_T(\un u)}))-1)\tau_{N_T(\un u)+1}} \label{density} \\
&=&
5^{N_T(\un{u})}
 \prod_{i=0}^{N_T(\un{u})-1} e^{-(h(\un{u}(t_i))-1)\tau_{i+1}} \lambda_{ \un{u}(t_i) } (\un{u}(t_{i+1}) - \un{u}(t_i))\times \nonumber\\
&&h(\un u(t_{N_T(\un u)}))e^{-(h(\un u(t_{N_T(\un u)}))-1)\tau_{N_T(\un u)+1}} \nonumber
\end{eqnarray}
where $\un{u}(\cdot) \in X_T$ with $N_T(\un{u})$ jump moments $0=t_0<t_1<\cdots<t_{N_T(\un{u})}<t_{N_T(\un{u})+1}=T$. For any  $\un u (\cdot) \notin X_T$,  ${\cal P} \left( \un{u}(\cdot) \right) = 0$.
Hence
\begin{equation}\label{dens}
  {\bf P} (\xi(\cdot)\in E)= e^T
  {\bf E} (e^{-A_T(\zeta)+B_T(\zeta)+ N_T(\zeta)\ln 5 };~
  \zeta(\cdot)\in E)
\end{equation}
for any measurable set $E\subseteq X_T$, where for $\un u \in E$
\begin{eqnarray}
 A_T (\un{u}) &:=& \sum_{i=0}^{N_T(\un u)} h(\un{u}(t_i))\tau_{i+1}=
    \int_0^Th(\un{u}(t))dt, \label{AT} \\
  B_T(\un{u}) &:=& \sum_{i=0}^{N_{T} (\un u)-1}\ln \bigl(\lambda_{ \un{u}(t_i) } ( \un{u}(t_{i+1}) - \un{u}(t_i)) \bigr)+\ln h\bigl(\un u(t_{N_T(\un u)})\bigr). \label{BT}
\end{eqnarray}

We study an asymptotic behavior of the logarithm of the probability ${\bf P} \bigl(\xi_T(\cdot)\in U_{\varepsilon}(\un{f}) \bigr)$ for any $\un f\in F$ using \reff{dens}.
The main contribution in this asymptotic comes from $A_T$. To prove this we consider
the scaled processes $\zeta_T(s)=\frac{\zeta(sT)}{T}, s\in [0,1]$.  Let $\un{x}(s)=\frac{\un{u}(sT)}{T}$ for $\un{u}\in X_T$, then
\begin{eqnarray*}
   A_T(\un{x}) := A_T(\un u) &=& T^2\int_0^1\left[\frac{c_0}{T}+c_1\frac{u_1(sT)}{T}+c_2\frac{u_2(sT)}{T}+c_3\min\Bigl\{ \frac{u_1(sT)}{T}, \frac{u_2(sT)}{T} \Bigr\}\right]\ed s \\
    &=&T^2\int_0^1\left[\frac{c_0}{T}+c_1x_1(s)+c_2x_2(s)+c_3\min\left\{x_1(s),x_2(s)\right\}\right]\ed s \\
    &=& T^2\left[\frac{c_0}{T}+I(\un{x})\right].
  \end{eqnarray*}
Then for any $\varepsilon$ there exists $\delta$ such that
\begin{equation}\label{3}
T^2 I(\un{f}) (1- \delta)  \le A_T(\un{x}) \le T^2I(\un{f}) (1+\delta)
\end{equation}
for any $\un{x}\in U_{\varepsilon}(\un{f})$. Hence
\begin{equation}\label{44.1}
L_+\le -I(\un{f})+\limsup_{\varepsilon\to 0}\limsup_{T\to \infty}\frac{1}{T^2}\ln {\bf
E} \bigl(e^{B_T(\zeta) + N_T(\zeta) \ln 5};~~\zeta_T(\cdot)\in  U_{\varepsilon}(\un{f}) \bigr).
 \end{equation}
Next we show that the second term in  (\ref{44.1}) is equal to 0.  

Let $\un{y}\in  U_{\varepsilon}(\un{f})$ 
and $K_+=K_+(\un{y})$ be the number of jumps of $\un{y}(\cdot)=(y_1(\cdot),y_2(\cdot))$  
on the time interval $\left[0,1\right]$,
such that the values of either $y_1$ or $y_2$ are increasing at the jump moments.
Recall that the path $\un{y}$ can increase by the increments $(1,0)$ or $(0,1)$.

Let $\varepsilon >0$ be such that $f_i(1) - \varepsilon >0, i=1,2$, then $y_i(1)>0$, since $\un y \in U_{\varepsilon}(\un f)$.
Thus
$$
K_+-K_- > 0
$$
where $K_-$ is the number of jumps on the time interval $[0,1]$, when the values of either $y_1$ or $y_2$ or both are decreasing at the jump moments. Note that $N_T(\un{y}) = K_+ + K_-$, and hence
\begin{equation}\label{7.2}
K_+ > \frac{1}{2}N_T(\un{y}).
\end{equation}
 
The next step of the proof is based on the following lemma.
 \begin{lemma}\label{Lema}
For any $\un f\in F$ there exist  positive constants $R_1$ and $R_2$, which depends on $\un{f}$, such that 
\begin{equation}\label{46.1}
e^{C_T} := {\bf E} \bigl(e^{B_T(\zeta) + N_T(\zeta) \ln 5};~\zeta_T(\cdot)\in  U_{\varepsilon}(\un{f}) \bigr) 
 \le  {\bf E} \exp \Bigl\{ 
\frac{N_T(\zeta)}{2} 
 (\ln T+R_1)+ \frac{1}{2} \ln(R_2T) \Bigr\}
 \end{equation}
 holds for small $\varepsilon$ (see \reff{44.1}).
 \end{lemma}
 \proof Let $\un{x}$ be some scaled trajectory of unscaled path $\un{u} \in X_T$, $\un{x}(s)=u(sT)/T, s\in [0,1]$ and $\{\tilde s_i\} \subset \{ s_i\} = \{ t_i / T \}$ be a subset of moments when the values $x_1$ or $x_2$ or both are decreasing. Remember that the number of such jumps is $K_-$. Thus  (see \reff{BT} for the definition of $B_T$):
\begin{eqnarray}
B_T(\un{u}) &:=& B_T(\un{x}) = \sum_{i=0}^{N_T(\un{x})-1} \ln \bigl(\lambda_{ T\un{x}(s_i) } \left( T(\un{x}(s_{i+1}) - \un{x}(s_i)) \right) \bigr)+ \ln \left( h( T \un{x}(t_{N_T(\un{x})}) ) \right) \nonumber\\
&\le & K_+ \ln c_0 + (K_-+1) \ln \Bigl( T \bigl( c_0 + (\max c_i) \bigl( \sup_{t\in [0,1]} \max\{ f_1(t), f_2(t)\} + \varepsilon \bigr) \bigr) \Bigr) \nonumber \\
& \le & \frac{1}{2}(N_T(\un{x})+1) ( \ln T + C ), \label{estBT2}
\end{eqnarray}
for some constant $C$ that depends on $\un{f}$. Choosing
$R_1 = C + 2\ln 5,~R_2=e^C$ we obtain the proof of the lemma. \qed 

To finish the proof of
$$
  \limsup_{\varepsilon \to 0} \limsup_{T\to \infty}
  \frac{1}{T^2}
  C_T  = 0,
$$
remark that the random variable $N_S(\zeta)$ has
Poisson distribution with a parameter $S$. Hence
$$
  {\bf E}e^{\theta N_S(\zeta)}=e^{S(e^\theta-1)}.
$$
Using \reff{46.1} we obtain
$$ 
e^{C_T}\le
e^{T(e^{\frac{1}{2}(\ln T+R_1)}-1)}R_2T\le
e^{T^{3/2}e^{\frac{R_1}{2}}}R_2T,
$$
which implies that
\begin{equation}\label{oeps}
  \limsup_{T\to \infty}
  \frac{1}{T^2}
  C_T\le \lim_{T\to\infty}\frac{1}{T^2}\left(
 T^{3/2}e^{\frac{R_1}{2}} + \ln (R_2T) \right)= 0.
\end{equation}
Therefore the proof of the upper bound \reff{2} is completed.

\vspace{1cm}

\noindent \underline{Lower Bound}. We have to prove the inequality
 \begin{equation}\label{4}
L_-:=\liminf_{\varepsilon\to 0} \liminf_{T\to \infty}\frac{1}{T^2}\ln
  {\bf P} \bigl( \xi_T\in  U_{\varepsilon}(\un{f}) \bigr)\ge -I(\un{f}).
\end{equation}
The probability of the event ${\mathcal U}(\un f) := (\xi_T\in  U_{\varepsilon}(\un{f}))$ can be lower estimated
by the probability of a more restricted event ${\mathcal U} (\un f, C) := (\xi_T\in U_{\varepsilon}(\un{f}),~N_T(\xi)\le CT)$.
A value of the constant $C$ depends on $\un{f}$. 
Using the
representation of the distribution of $\xi$ in terms of the process $\zeta$ (see \reff{dens}), the inequalities \reff{3} and that $B_T(\un x)> N_T(\un x) \ln(\tilde c)$, where $\tilde c :=\min c_i$, we obtain the
lower estimate
\begin{eqnarray}\label{4.1}
 &&\liminf_{T\to \infty}\frac{1}{T^2}\ln
  {\bf P} \bigl( \xi_T\in U_{\varepsilon}(\un{f}) \bigr) \\
 && {} \ge -I(\un{f})(1+\delta)+
  \liminf_{T\to \infty}\frac{1}{T^2}\ln
 {\bf E} \bigl( e^{N_T(\zeta) \ln(5 \tilde c)}; \zeta_T\in U_{\varepsilon}(\un{f}),\, N_T(\zeta)\le CT \bigr).\nonumber
\end{eqnarray}
If $\ln(5 \tilde c) >0$, then $e^{N_T(\zeta) \ln(5\tilde c)} >1$ and the expectation in \reff{4.1} is bounded below by the probability ${\bf P} ({\mathcal U} (\un f, C) )$. On the other hand, if  $\ln(5 \tilde c) <0$, then the expectation is bounded below by $e^{CT \ln(5 \tilde c)} {\bf P}({\mathcal U} (\un f, C))$.

Recall that on the event $\mathcal U (\un f, C)$, the values of the process $\zeta_T(t)$ are non-negative. The lower estimate  of $\ln{\bf P}({\mathcal U} (\un f, C))$ follow from the recent result in \cite{BM1}, (Theorems 3.1 and 3.3). Namely, there exists a constant $J>0$ such that
\begin{equation*}
\liminf_{T\to\infty} \frac{1}{T}\ln {\bf P} \bigl( \zeta_T\in  U_{\varepsilon}(\un{f}),\, N_T(\zeta_T)\le CT \bigr) \geq  J> -\infty.
\end{equation*}
Thereby
\begin{equation*}
\liminf_{T\to \infty}\frac{1}{T^2}
\ln {\bf P} \bigl(\zeta_T\in U_{\varepsilon}(\un{f}),\, N_T(\zeta_T)\le CT \bigr)=0.
\end{equation*}

Despite that the formula for the rate function \reff{rate} can be applied for discontinuous functions, in fact the rate function is infinite 
for such functions. That happens because 
$ {\bf P} (\xi_T \in U_\varepsilon( \un x) ) = 0$ in the
uniform topology   for any discontinuous function $\un x$ if $\varepsilon$ is small enough. 
\qed

\begin{remark}
In \cite{BM1}, the large deviation principle is proved for real valued processes with independent increments. The result of Theorems 3.1 and 3.3 from \cite{BM1} can be easily extended to finite dimension cases. 
\end{remark}

\subsection{A version for ``integral" large deviation principle.}

For any continuous function $\un{f}=(f_1,f_2)\in F$ and any positive $\varepsilon$ and $M$, consider the following sets:
\begin{eqnarray*}
B_{\un{f},\varepsilon,M} &=& \{ \un{x}=(x_1,x_2) \in X: \  f_i (t)-\varepsilon \le x_i(t) \le M, \ i=1,2, \ t\in [0,1] \}, \\
B_{\un{f},M} &=& \{ \un{x}=(x_1,x_2) \in X: \  f_i (t) \le x_i(t) \le M, \ i=1,2, \ t\in [0,1] \}.
\end{eqnarray*}
We will call them strips.

\begin{theorem}\label{thstrip} For any $\un{f} \in F$ and any $M> \sup_{t\in [0,1]} \max \{ f_1(t), f_2(t) \}$
\begin{equation}\label{ldRed}
\lim_{\varepsilon \to 0}  \lim_{T\to \infty}\frac{1}{T^2}\ln {\bf P} \left( \xi_T(\cdot) \in B_{\un{f},\varepsilon,M}  \right) =
- \inf_{\un{g}\in F \cap B_{\un{f},M} } I(\un{g} )= -I( \un{f} ).
\end{equation}
\end{theorem}
\proof  The upper bound follows from representation \reff{dens}: for any $\varepsilon$ there exists $\delta$ such that
\begin{equation*}
\frac{1}{T^2}\ln {\bf P} (\xi_T(\cdot)\in B_{\un{f},\varepsilon,M} ) \leq -\inf_{\un{g}\in F \cap B_{\un{f},M} }I(\un{g})(1-\delta) + \sup_{\un{g}\in F \cap B_{\un{f},M} } \frac{1}{T^2} C_T,
\end{equation*}
see Lemma~\ref{Lema} for the definition of $C_T$. The proof of the relation
$$
\limsup_{T\to\infty} \sup_{\un{g}\in F \cap B_{\un{f},M} } \frac{1}{T^2} C_T = o(1) \mbox{ as }\varepsilon \to 0
$$
basically repeats the arguments of Section \ref{22} replacing $\sup_{t\in [0,1]}  \max \{ f_1(t), f_2(t) \}$ by $M$ in  \reff{estBT2}.  This modification does not affect the principal inequality \reff{oeps}. It proves the upper bound.

The lower bound becomes obvious using the estimation 
$$
{\bf P} \bigl( \xi_T(\cdot) \in B_{\un{f},\varepsilon,M}  \bigr) \ge {\bf P} \bigl( \xi_T(\cdot) \in U_{\varepsilon}(\un{f})  \bigr),
$$
and after that the usage of Theorem~\ref{thlldp} completes the proof of the theorem. \qed

\begin{remark}
Theorem~\ref{thstrip} holds also if, in the definition of the strip, we substitute the upper bound $M$ by a bound  $(M_1, M_2)+\un{g}$, where $\un{g}$ is any continuous function on [0,1] with $\un{g}(0)=(0,0)$ and $M_1, M_2$ are some positive constants. The lower bound is defined by a function $\un{f} \in F$ such that
$$
\sup_{t\in [0,1]} \bigl( M_i + g_i(t) - f_i(t) \bigr) >0.
$$
for any $i=1,2.$
\end{remark}

We have not proven the large deviation principle in its complete form. There are some reasons for this.
First, the rate function \reff{rate} is not compact. Second, in the considered topology the exponential tightness does not hold. Moreover, the space $X$ is not complete and it is not separable. Thus we stated the large deviation for some special sets, that we  called the strips. It seems, the strips can be demanded in applications.

\section{Acknowledgments} The authors thank N. Vvedenskaya for a number of useful discussions, 
and A.A. Borovkov for stimulating questions.

The work of A.M. was partially supported by grant FAPESP (2012/07845-3), grant of President of RF (NSh-3695.2008.1) and of RFFI (08-01-00962). The work of E.P. was partially supported by grant of RFFI (11-01-00485). A.Y. thanks The National Council for Scientific and Technological
Development (CNPq), Brazil, grant 308510/2010-0. E.P. thanks University of S\~{a}o Paulo (USP) and NUMEC for warm hospitality.

\end{document}